\documentclass[12pt,letterpaper, reqno]{amsart}

\usepackage{amssymb,latexsym, amsmath, amsxtra}
\usepackage[dvips]{graphics}
\usepackage{hyperref}

\theoremstyle{plain}
\newtheorem{theorem}{Theorem}[section]

\theoremstyle{definition}
\newtheorem{definition}[theorem]{Definition}
\theoremstyle{remark}

\numberwithin{equation}{section}
\numberwithin{theorem}{section}
\numberwithin{table}{section}
\numberwithin{figure}{section}

\newcommand{\conj}[1]{\overline{#1}}
\newcommand{\term}[1]{\emph{\textbf{#1}}}

\def\({\left(}
\def\){\right)}

\begin{document}
\title{When are the roots of a polynomial real and distinct? \\ A graphical view}
\author{David~W.~Farmer}
\address{American Institute of Mathematics,
600 East Brokaw Road,
San Jose, CA 95112-1006}
\email{farmer@aimath.org}


\begin{abstract}
We prove the classical result, which goes back at least to Fourier,
that a polynomial with real coefficients has all zeros
real and distinct if and only if the polynomial and also all of its
nonconstant derivatives have only negative minima and 
positive maxima.  Intuition for the result, involving illuminating pictures,
is described in detail.
The generalization of Fourier's theorem to certain entire functions of
order one (which is conjectural) suggests that
the official description of the Riemann Hypothesis Millennium Problem
incorrectly describes an equivalence to the Riemann Hypothesis.
The paper is reasonably self-contained and is intended be accessible
(possibly with some help) to
students who have taken two semesters of calculus.
\end{abstract}

\keywords{real polynomial, real roots, distinct roots, hyperbolic polynomial, Riemann Hypothesis}
\thanks{This research was supported by the National Science Foundation.}

\maketitle

\section{Introduction} 

A recent paper in the Monthly~\cite{Cha}
\emph{When are the roots of a polynomial real and distinct?} provided
the following elegant criterion:
\begin{theorem}\label{thm:inequality}
Let $f$ be a polynomial of degree $N\ge 1$ with real coefficients.
Then the zeros of $f$ are real and distinct if and only if
\begin{equation}\label{eqn:maininequality}
(f^{(j)}(x))^2 - f^{(j-1)}(x) f^{(j+1)}(x) > 0
\end{equation}
for all $x\in \mathbb R$ and all $1\le j \le n-1$.
\end{theorem}
The notation $f^{(j)}$ refers to the $j$th derivative of~$f$.

Suppose one had a polynomial and wanted to know if its roots are
real and distinct.  
How could one actually check the condition in Theorem~\ref{thm:inequality}?
One possibility is to graph \eqref{eqn:maininequality} for each~$j$,
and observe whether the graphs always stay above the $x$-axis.
That idea was the motivation for seeking a condition that involved graphing the
polynomial itself.

\begin{figure}[htp]
\scalebox{0.6}[0.6]{\includegraphics{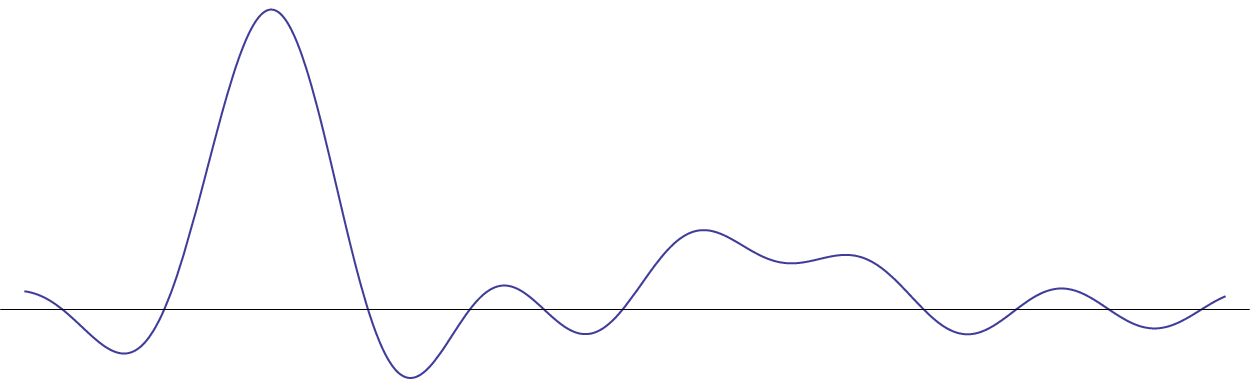}}
\hskip 0.4in
\scalebox{0.6}[0.6]{\includegraphics{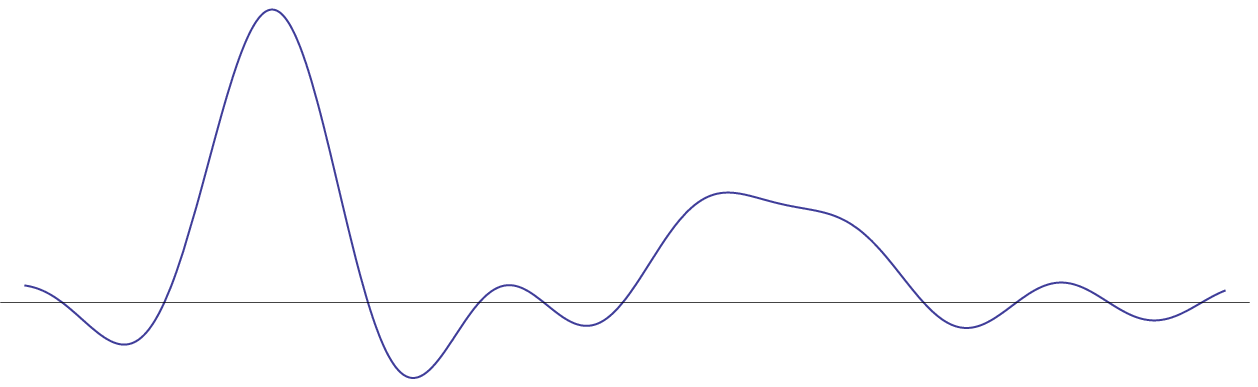}}
\vskip 0.4in
\scalebox{0.6}[0.6]{\includegraphics{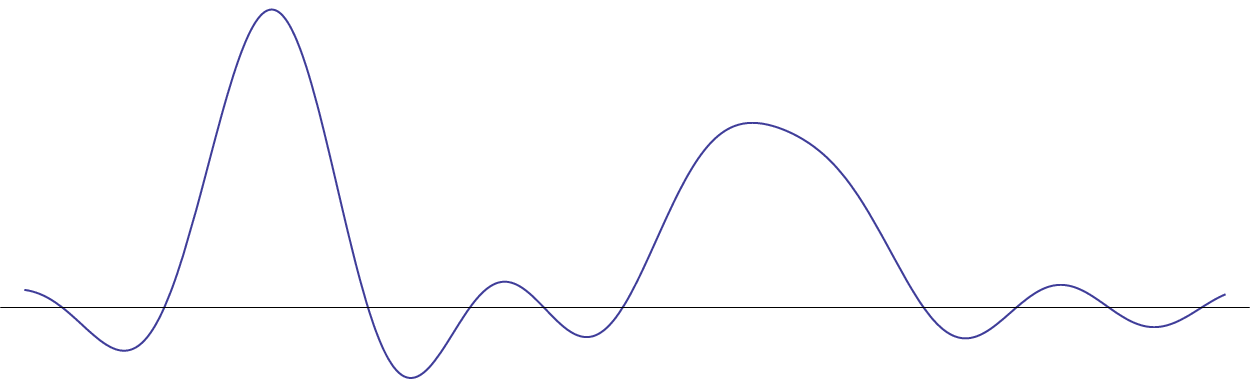}}
\hskip 0.4in
\scalebox{0.6}[0.6]{\includegraphics{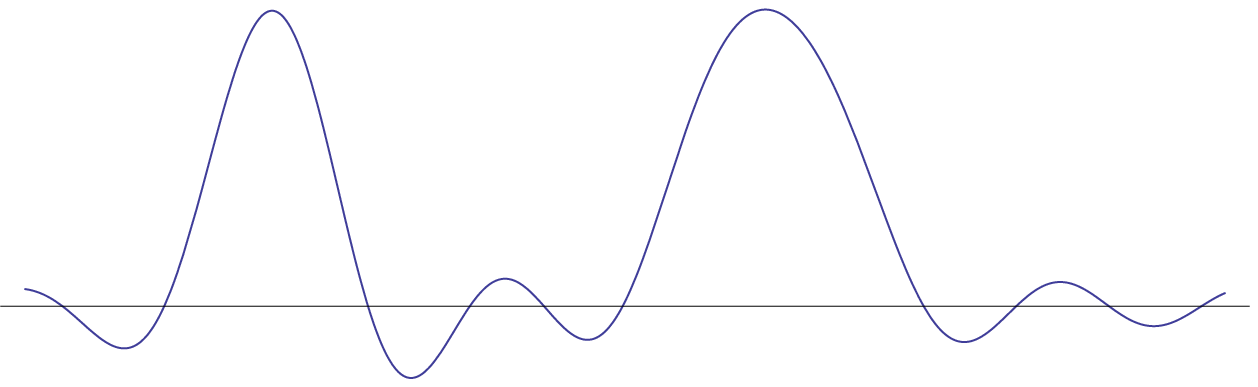}}
\caption{\sf
Graphs of portions of four different high-degree polynomials,
each of which has a pair of non-real zeros in the given range.
} \label{fig:motivatingplots}
\end{figure}

To help build intuition,
consider the plots in Figure~\ref{fig:motivatingplots},
Each is a graph of part of a different high-degree polynomial which does \emph{not} have
only real zeros.  Try to decide if the information in each graph is sufficient
to conclude that the polynomial has a non-real zero.
Note that the graphs show the $x$-axis, making it possible to see
the zeros of the polynomials, but both the horizontal and vertical
scales are omitted because those are irrelevant to the discussion.

We review some basic facts about polynomials, and then discuss the
clues which can be seen in the graphs in Figure~\ref{fig:motivatingplots}.

\section{Algebraic and geometric properties of polynomials}

\term{Polynomials} are expressions of the form $f(x) = a_0 + a_1 x + \cdots + a_N x^N$.
The $a_j$, which are (real or complex) numbers, are the \term{coefficients}
of the polynomial.  If $a_N \not= 0$, then it is called the 
\term{leading coefficient}, and we call $N$ the \term{degree}
of the polynomial. We write $\deg(f) = N$.
A \term{root}, also called a \term{zero}, of the polynomial
$f$ is a (real or complex) number $z$ such that $f(z)=0$.

If the degree of $f$ is at least one, then we say $f$ is a \term{nonconstant}
polynomial.  If $f$ is a nonconstant polynomial of degree~$N$ then
its derivative $f'$ is a polynomial of degree~$N-1$.  That is the first of
several statements we will encounter where a natural-looking assertion
is actually incorrect when applied to a constant polynomial.  Another
edge case is the \term{zero polynomial} $f(x) = 0$.  It is a polynomial,
but its degree is not defined.

The Fundamental Theorem of Algebra states that a degree $N$ polynomial
can be factored into a product of $N$ linear factors:
\begin{equation}\label{eqn:polyfactored}
f(x)=a_N(x-z_1)\cdots (x-z_N).
\end{equation}
The numbers $z_1$, \ldots, $z_N$ are the
roots of the polynomial.  If a root $z_j$ appears $k$ times then we say
that the root $z_j$ has \term{multiplicity}~$k$.  It follows from the
product rule that if a root of $f$ has multiplicity greater than~1,
then it is also a root of the derivative~$f'$, with multiplicity
decreased by one.
We observe that it takes $N+1$ numbers to specify a polynomial.
Those numbers could be the coefficients, or they could be the
leading coefficient and the zeros.

For the remainder of this paper, we will assume that all polynomials
have real coefficients.  We will refer to these as \term{real polynomials}.
A real polynomial defines a
differentiable function on the real numbers, so we can graph it,
as illustrated in Figure~\ref{fig:motivatingplots}.
The zeros of a real polynomial are either real, or they occur in
complex conjugate pairs:  if $z = x + i y$ is a root (here $x$ and $y$
are real, and $i$ is the imaginary unit, which satisfies $i^2 = -1$),
then its complex conjugate $\conj{z} = x - i y$ is also a root,
and those roots have the same multiplicity.  We say that the zeros of
a polynomial are \term{distinct} if every root has multiplicity one.
If we refer to ``the distinct real zeros'' of $f$, we mean the
list of roots where each root of higher multiplicity only appears once.
If we wish to refer to all zeros and want to be particularly clear that
we really mean \emph{all} zeros, we may say ``zeros counted with multiplicity.''

Since real polynomials are continuous and have continuous derivatives,
Rolle's theorem implies that between every consecutive pair of
real zeros, $(z_j, z_{j+1})$, with $z_j \not= z_{j+1}$,
there is a zero of the derivative.  
Note that if $f$ has a multiple zero, meaning $z_j = z_{j+1}$,
the conclusion of Rolle's theorem still holds,
because of our previous observation on differentiation and multiple zeros.
However, we will not consider the case of multiple zeros in this paper.

Since the polynomial does not
change sign on the interval $(z_j, z_{j+1})$, on that interval either
the polynomial is positive and it has a local maximum,
or it is negative and it has a local minimum.
(It may also have other critical points on that interval.)
In particular, if the
polynomial $f$ has $n$ distinct real zeros, then its derivative $f'$ has at least
$n-1$ distinct real zeros.  It will have exactly $n-1$ distinct real zeros
if and only if there is
exactly one critical point between each adjacent pair of real zeros,
and no critical points for $x < z_0$ or $x > z_N$.  
In that situation $f'$ could have $n-1$ distinct real zeros but more than $n-1$ real zeros
counted with multiplicity.
That would happen, for example, if there was a local maximum where the 2nd and 3rd derivatives were also zero.

Since $f'$ has one fewer
zero than~$f$, the above discussion proves:

\begin{theorem}\label{thm:easydirection} If a real nonconstant polynomial has only distinct real zeros,
then the same is true of its derivative.  Furthermore, for each nonconstant
derivative, every local maximum is positive, and every local minimum is negative.
\end{theorem}

The plots in Figure~\ref{fig:motivatingplots} are meant to suggest
the converse of Theorem~\ref{thm:easydirection}:

\begin{theorem}\label{thm:harddirection} Suppose $f$ is a real nonconstant polynomial
such that it and all of its 
nonconstant derivatives have  every local maximum is positive, and every local minimum is negative.
Then $f$ has only real zeros, and all zeros are distinct.
\end{theorem}

Theorem~\ref{thm:harddirection} goes back at least to
Joseph Fourier (1768-1830).  We give Fourier's proof in Section~\ref{sec:fourier}.

\section{Intuition for non-real roots}
We build some intuition for Theorem~\ref{thm:harddirection} before proving it
by induction in the next section,
and then giving a more conceptual proof in  Section~\ref{sec:fourier}.

The main idea is that if $f$ does not have all real zeros, the non-real
zeros cause \term{funny business} in the graph of~$f$.  A good eye for
detecting funny business requires looking at plots of many polynomials.
The examples in Figure~\ref{fig:motivatingplots} were chosen to
illustrate the main ideas.

Let's start with the upper-left plot.  One can count that there are 10 real
zeros in the given domain.  The graph goes up and down in a pleasant
way, constrained by the fact that it must cross the axis at the real zeros.
Polynomials don't like to change direction too quickly, so when there
is a big gap between zeros, the polynomial moves further from the $x$-axis
before turning around.  The gap between adjacent real zeros is the most
important factor in determining how far the graph wanders from the axis,
but it is not the only factor.  For example, the gap between the 1st and 2nd
zeros is the same as the gap between the 3rd and the 4th, yet the graph
becomes more negative between that second pair.  This is due to the
influence of the neighboring zeros.  The minimum being less negative between
the 1st and 2nd zero is a clue
about the zeros which are slightly to the left of the displayed range.    

But what is going on between the 6th and 7th real zeros?
Why does the graph start wiggling for no good reason?  There is a good
reason for that wiggling:  it is caused by a pair of conjugate non-real zeros of the polynomial.
That behavior is the most blatant example of funny business:  a positive local
minimum.  From the graph we can deduce some information about those
non-real zeros:  its real part is slightly closer to the 7th zero than
to the 6th zero, and its imaginary part is not too large (on the scale
of the gap between the adjacent real zeros).

Let's move on to the upper-right plot.  Again we see funny business
between the 6th and 7th real zeros, although it is not as extreme
as in the previous example.  It is also caused by a pair of conjugate
zeros, but those zeros have slightly larger imaginary part
than in the previous graph
(in fact, 50\% larger).  There is only one critical point in that
interval, in particular no positive local minimum.  But what will happen
when we take the derivative of that function?  Looking carefully at
the graph, we can see that the derivative will have a negative local
maximum.

Now consider the lower-left plot.  The only thing which has changed is
the imaginary part of the complex roots, increasing by another 30\%.
Again we see funny business, in the form of an odd bulge to the graph
as it moves from the local maximum down to the 7th real zero.
That type of funny business might not have an official name,
but if you are alert to look out for it, you can immediately
identify it as coming from non-real zeros.  If you differentiate,
the funny business becomes more pronounced.  The second derivative
has a positive local minimum.

The reason the funny business becomes more extreme with each derivative
(until all the zeros end up on the real axis and the funny business
stops) is that the imaginary
part of the non-real zeros decreases with each derivative.  Zeros
with large imaginary part cause a slight disturbance, while zeros
with small imaginary part have a greater effect.  
The size of the imaginary part causing funny business
should be measured relative to the
spacing between the nearby zeros.

Finally, let's look at the lower-right graph. The only change is that
the imaginary part of the complex roots increased by another 30\%.
Can we see any funny
business there?  There might be some asymmetry between the 6th and 7th
zeros, but maybe not enough to draw a definite conclusion.  
There is even more asymmetry between the 1st and 2nd zeros, 
and that has nothing to do with non-real zeros.
Looking more closely there are other clues.  The height of the
maximum between the 6th and 7th zeros is about the same as
between the 2nd and 3rd, but the gap between the 6th and 7th is
much larger.  That contradicts our previous observation about
the gap being the dominant factor in determining the size of the
maximum.  So, something (i.e., a conjugate pair of non-real zeros)
is causing the maximum to be smaller than expected.

There is another clue.  Look closely at the shape of the graph near the
two largest local maxima.  The one on the left is more sharply curved, and more
quickly becomes approximately straight.  The one on the left is what
the graph of a polynomial with only real zeros should look like.
The one on the right is an example of funny business.
The 3rd derivative of that polynomial has a positive local minimum.

%
%

\section{Induction proof of Theorem~\ref{thm:harddirection}}

The proof will proceed by induction on the degree of the polynomial.
A different and more interesting proof,
using ideas of Fourier,
is in Section~\ref{sec:fourier}

Suppose $f$ is a degree 1 polynomial.  It has no non-constant derivatives,
so the conditions in the statement are vacuously true, so we want to conclude
that $f$ has only real zeros.  That is true because the graph of $f$ is a line
of nonzero slope (nonzero because of the definition of `degree') so the
graph must cross the $x$-axis.  Thus, $f$ has at least one real root, 
which is all roots by the Fundamental Theorem of Algebra,
and the root is distinct because there is only one of them.

For the induction step we
will prove the contrapositive, which is that if $f$ does not have
all distinct real roots, then either $f$ or some derivative of $f$ must have
a local maximum that is non-positive, or a local minimum
that is non-negative.  
So, suppose the theorem is true for all real polynomials of degree~$N$,
and suppose $f$ is a real polynomial of degree~$N+1$ which either has
a non-real root, or it has a real root that is non-distinct.

In the case of a non-distinct real root,  if the multiplicity is even,
then the root is either a local maximum or local minimum.  Since the
function is zero at that root, the maximum or minimum is non-positive
and non-negative, which is what we wanted to show. 
(This situation will appear repeatedly, so we will abbreviate it
as ``has a zero maximum or minimum''.) If the multiplicity
is odd, the multiplicity must be at least 3, so the derivative
has a multiple root of even multiplicity, which again is what we wanted
to show.

Now we handle the case that the real roots of $f$ are distinct,
but not all roots are real.  If $z_1 \le z_2 \le \cdots \le z_M$ are the real roots of~$f$,
then $M\le N-1$ because the non-real roots of $f$ come in conjugate
pairs, so there must be at least two of them.

Suppose $f'$ has $M-1$ distinct real roots, which it the minimum
allowed by Rolle's theorem.  If those roots all have multiplicity one,
then $f'$ has non-real roots because the degree of $f'$ is at least~$M+1$.
Since the degree of $f'$ is $N$, this case is covered by the induction hypothesis.

Now suppose $f'$ has $M-1$ distinct real roots, at least one of which
has multiplicity greater than one.  Then as we argued previously,
either $f'$ or $f''$ has a zero maximum or minimum.

This leaves the final case, where $f'$ has at least $M$ distinct real zeros.
Either there must be a zero of $f'(x)$ outside the interval $(z_1, z_M)$,
or by the pigeon hole principle
at least two zeros in an interval $(z_{j}, z_{j+1})$.

In the case of a zero outside $(z_1, z_M)$, by replacing $f(x)$ by
$f(-x)$ and/or by $-f(x)$ we may assume the zero of $f'(x)$ lies in $x > z_M$,
and furthermore $f(x) > 0$ for $x > z_M$.
If that zero of $f'$ was a local minimum of $f$, we are done.  If it 
is a local maximum, then since $f$ is increasing for $x$ sufficiently large,
there must be a local minimum to the right of that local maximum,
so again we are done.  The final possibility is that it is neither,
in which case it corresponds to a multiple zero of $f'$,  so either
$f'$ or $f''$ has a zero maximum or minimum,
as in previous cases.

The final possibility is that $f'$ has at least two zeros in an interval $(z_{j}, z_{j+1})$.
If the two zeros of $f'$ are not distinct, then either $f'$ or~$f''$ has a
zero minimum or maximum, and we are done.  Thus,  we may assume those two zeros of $f'$ 
in that interval are distinct.
As in the previous case, we may assume $f$ is positive on that interval,
and so it has a local maximum, which corresponds to one of those
distinct zeros of $f'$.  If the other zero of $f'$ corresponds to a minimum of $f$,
we are done.  If it corresponds to a maximum, we are done because there must be a
minimum between two maxima.  That leaves the third possibility that 
the other zero of $f'$ corresponds to neither a maximum nor a minimum of~$f$,
but that means it is a multiple zero of $f'$, so again we are done.

This covers all cases in the induction step, so the proof is complete.

\section{Fourier's proof}\label{sec:fourier}
We now give a different proof, using ideas due to Joseph Fourier
of Fourier series fame.  The proof is logically the same, in the sense
that there is a direct correspondence between the ideas of Fourier
and the cases in the induction proof.  But 
Fourier's proof has better explanatory power, and it also illustrates
the value of introducing useful terminology.
This section is written as if we were reading Fourier's mind; the author
was introduced to these ideas by reading a paper of Y.-O.~Kim~\cite{Kim}.

\begin{figure}[htp]
\scalebox{1.3}[1.3]{\includegraphics{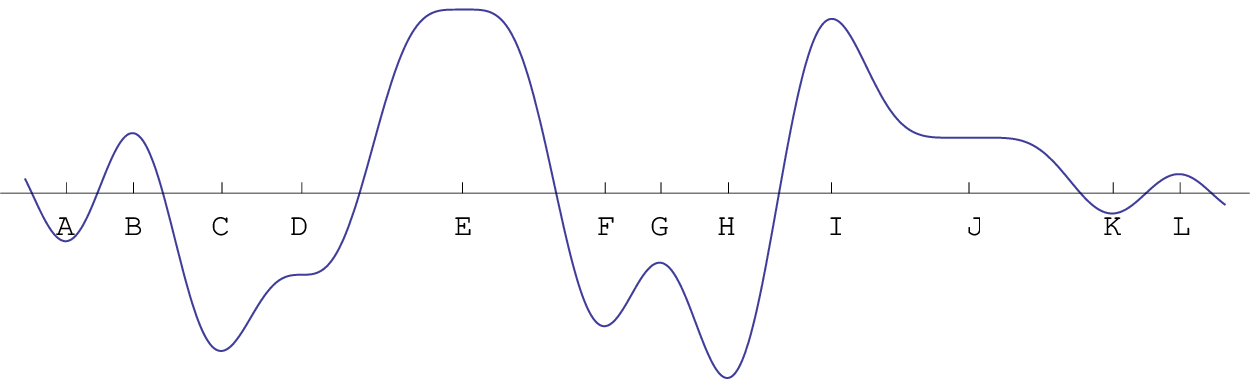}}
\vskip 0.4in
\scalebox{1.3}[1.3]{\includegraphics{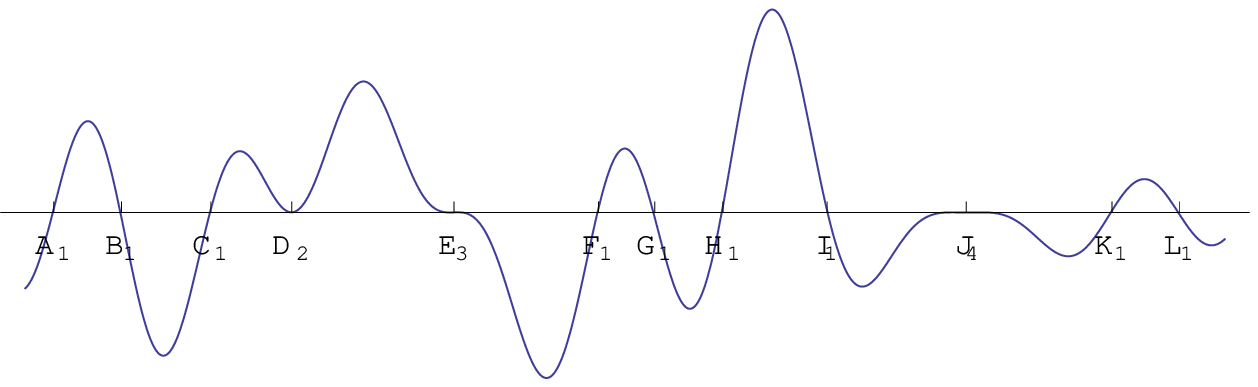}}
\caption{\sf
Part of the graph of a high degree polynomial (top),
and its derivative (bottom).
The labeled points $A$, $B$, \ldots , $L$ are critical points in the top graph, equivalently, zeros of the derivative in the bottom graph.
In the bottom graph the subscript indicates the multiplicity of the zero.
} \label{fig:fourier}
\end{figure}

Rolle's theorem for a polynomial $f$ implies that there is a zero of $f'$
between each consecutive pair of real zeros, but sometimes you get ``extra''
zeros of $f'$.  Figure~\ref{fig:fourier} illustrates some of the possibilities.

Fourier did not think all critical point were equally interesting.  For example,
critical point $B$ is not particularly interesting:  by Rolle's theorem
there has to be a critical point between the zeros on either side of $B$,
so there is nothing special about the point $B$.  Similarly,
none of $A$, $C$, $F$, $H$, $I$, $K$, or $L$ is interesting to Fourier.
The remaining critical points are interesting:  $D$, $E$, and $J$ because in
addition to $f'$ vanishing, so do some additional derivatives, so those
points contribute ``extra'' zeros to $f$; and point $G$ is interesting because
it is a negative local maximum, which is not predicted by
Rolle's theorem.  Table~\ref{tab:fourier} compiles information about the
behavior at each critical point.  The entries $+$, $0$, and $-$ refer to
the function being positive, negative, or zero; a $*$ means a nonzero
value (which could be positive or negative; and a blank entry means
the value is irrelevant in terms of Fourier's classification of critical points.

\begin{table}[]
\begin{tabular}{cccccccc}
$x$ &  $f(x)$  & $f'(x)$  & $f''(x)$ & $f'''(x)$ & $f^{(4)}(x)$ & $f^{(5)}(x)$ & $k_{f'}(x)$ \\
\hline
$A$   & $-$   & 0  &  $+$  &  & &&0 \\
$B$   & $+$   & 0  &  $-$  &  & &&0 \\
$C$   & $-$   & 0  &  $+$  &  & &&0 \\
$D$   & $*$   & 0  &  0  & $*$  & &&1  \\
$E$   & $+$   & 0  &  0  & 0  & $-$ &&1 \\
$F$   & $-$   & 0  &  $+$  &  &&&0  \\
$G$   & $-$   & 0  &  $-$  &  & &&1 \\
$H$   & $-$   & 0  &  $+$  &  &&&0  \\
$I$   & $+$   & 0  &  $-$  &  &&&0  \\
$J$   & $*$   & 0  &  0  & 0  & 0 & $*$ & 2\\
$K$   & $-$   & 0  &  $+$  &  &&&0  \\
$L$   & $+$   & 0  &  $-$  &  & &&0 \\
\end{tabular}
\caption{Signs of the function in Figure~\ref{fig:fourier} and
its derivatives.  The final column indicates the number of pairs of
non-real zeros which cause the behavior of the other entries in the row.
}\label{tab:fourier}
\end{table}

The final column of Table~\ref{tab:fourier} is a
number that expresses how ``interesting'' the critical point $x$ is,
in Fourier's view.  
More specifically, $2k_{f'}(x)$ is the number of ``extra'' zeros that $f'(x)$ has
arising from the point~$x$,
where by \emph{extra} we mean ``more than one would expect only
from Rolle's theorem.''  
The extra zeros might be at~$x$, which happens at $D$, $E$, and $J$.
Or they might be nearby: the second extra zero from $G$ is either $F$ or $H$,
although it is not meaningful to specify which one it is.

 Looking at the plots in Figure~\ref{fig:fourier},
we see that $f$ has 9 real zeros in that interval, so by Rolle's
theorem $f'$ should have at least 8 real zeros in that interval.
We can count that $f'$ actually has 18 real zeros counting multiplicity
(because $D$ has multiplicity~2,
$E$ has multiplicity~3, and $J$ has multiplicity~4).  There are 10
extra zeros, which is twice the sum of the values of $k_{f'}$ in
Table~\ref{tab:fourier}.

Another interpretation of $k_{f'}(\alpha)$ is that it is the number of pairs of complex
conjugate zeros of $f$ which are needed to explain the funny business at
the critical point~$\alpha$.

To relate to our earlier discussions, if $k_{f'}(\alpha) > 0$,
then either $f$ or $f'$ or $f''$
has a non-negative minimum or non-positive maximum at~$\alpha$.

Here is a formula for~$k_{f'}$.  Suppose $\alpha$ is a critical point
which is not a zero of $f$, and furthermore suppose
\begin{equation}
f'(\alpha) = f''(\alpha) = \cdots = f^{(m)}(\alpha) = 0
\end{equation}
and $f^{(m+1)}(\alpha) \not = 0$.  That is, exactly the first $m$
derivatives of $f$ vanish at~$\alpha$.  Then
\begin{equation}
k_{f'}(\alpha) = \begin{cases}
\frac{m}{2} & \text{if $m$ is even} \cr
\frac{m-1}{2} & \text{if $m$ is odd and $f(\alpha)f^{(m+1)}(\alpha) < 0$} \cr
\frac{m+1}{2} & \text{if $m$ is odd and $f(\alpha)f^{(m+1)}(\alpha) > 0$}. 
\end{cases}
\end{equation}
That formula agrees with the values in Table~\ref{tab:fourier},
and we leave it as an exercise to check the other cases.
We set $k_{f'}(\alpha) = 0$ if $f'$ is identically zero, 
or if $\alpha$ is not a zero of $f'$,
so $k_{f'}(\alpha)$ is defined for all $\alpha \in \mathbb C$,
but it is $0$ for almost all~$\alpha$.
We will call $k_{f'}(\alpha)$ the \emph{Fourier multiplicity of the
critical zero $\alpha$ of $f'$}.

Now we assemble the above ideas into Fourier's proof of
Theorem~\ref{thm:harddirection}, introducing two new functions
which help us keep track of the quantities of interest.  First is the
function $K(g)$ which adds up all of the Fourier multiplicities of $g$:
\begin{equation}
K(g) =  \sum_{\alpha \in \mathbb C} k_{g} (\alpha).
\end{equation}
The sum is over all complex numbers~$\alpha$, but recall that
$k_{g} (\alpha) = 0$ for all but finitely many~$\alpha$,
so really it is a finite sum.

We use $K(g)$ to express an observation we made earlier.
Since $k_{f'}(\alpha)$ counts the number of pairs of non-real
zeros of $f'$ which are causing the funny business near $x=\alpha$
in the graph of $f(x)$, we have:
\begin{equation}\label{eqn:countfunny}
\text{(number of non-real zeros of $f$)}
- \text{(number of non-real zeros of $f'$)}
= 2 K(f').
\end{equation}
The factor of $2$ of the right of~\eqref{eqn:countfunny} comes from the
fact that $k_{f'}(x)$ counts \emph{pairs} of non-real zeros.
As a formula, $Z_C(f) - Z_C(f') = 2 K(f')$, where
$Z_C(g)$ counts the number of non-real zeros of the nonzero polynomial~$g$
(if $g$ is identically zero, we set $Z_C(g)=0$).
That was Fourier's key insight into the proof of Theorem~\ref{thm:harddirection}.
Rearranging and then using induction, we have
\begin{align}
Z_C(f) =\mathstrut& Z_C(f') + 2 K(f') \cr
Z_C(f) =\mathstrut& Z_C(f'') + 2 K(f') + 2 K(f'')\cr
\vdots\ \,& \cr
Z_C(f) =\mathstrut& Z_C(f^{(n)}) + 2 \sum_{j=1}^n K(f^{(j)}) 
\end{align}
for any $n \ge 1$.

If $n$ is larger than the degree of $f$ then $f^{(n)}$ is identically zero
so $Z_C(f^{(n)})=0$ and we have $Z_C(f) = 2 \sum_{j=1}^n K(f^{(j)})$.
If $f$ does not have all real zeros, then $Z_C(f) > 0$, so 
$ K(f^{(j)}) > 0$ for some $j$.  But that implies 
at least one of $f^{(j-1)}$, $f^{(j)}$, or $f^{(j+1)}$ has a non-negative
minimum or non-positive local maximum.  That completes Fourier's
proof of Theorem~\ref{thm:harddirection}.

The reader is invited to locate where in Fourier's proof we used
the fact that the roots of $f$ are distinct.

\section{Infinitely many zeros}

A question on MathOverflow~\cite{MO} asked
about an equivalence to the Riemann Hypothesis stated in the
official description~\cite{Cl} of the Riemann Hypothesis Millennium Problem:
\begin{quote}
\emph{
The Riemann hypothesis is equivalent to the statement that all local maxima
of $\Xi(t)$ are positive and all local minima are negative, \ldots.
}
\end{quote}
Here $\Xi(t)$, called the Riemann $\Xi$-function,
is closely related to the Riemann zeta-function.
See Section~\ref{sec:zeta} for definitions and background.
For this immediate discussion, it is sufficient to know that
$\Xi(t)$ is real for real~$t$, and the Riemann Hypothesis
is the assertion that all zeros of $\Xi(t)$ are real.
So, it is similar to the real polynomials we have been discussing,
except it has infinitely many zeros.  

The condition for a polynomial to have only real distinct zeros
involved maxima and minima of the original function, and also
of all derivatives; however, the quote above does not
mention derivatives.  Is our understanding of polynomials giving
us good intuition for functions with infinitely many zeros,
and so we should be skeptical of that claimed equivalence to
the Riemann Hypothesis?

As we will explain,
our intuition for polynomials is very likely guiding us in the
right direction.  However, establishing an analogue of Theorem~\ref{thm:harddirection}
for a class of functions which includes the Riemann $\Xi$-function,
is an unsolved problem!

Understanding the similarities between polynomials and certain
(but not all)
functions with infinitely many zeros,
and the relevance to the Riemann Hypothesis,
is the topic of the remainder of this paper.

\subsection{The Riemann zeta-function and the prime numbers}
The Riemann $\zeta$-function is defined by the Dirichlet series
\begin{equation}\label{eqn:zetaseries}
\zeta(s) = \sum_{n=1}^\infty \frac{1}{n^s} .
\end{equation}

That series was first studied by Euler, who used it to prove that there are
infinitely many prime numbers.  The infinitude of primes had already been proven by Euclid,
but Euler's result is
better:  he showed that the series $\sum_{n=1}^\infty 1/p_n$ diverges.
Here $p_1, p_2, p_3,\ldots = 2, 3, 5, \ldots$ is the sequence of prime numbers.
If the sequence of primes increased very quickly (a possibility that is
allowed by Euclid's proof), then Euler's series $\sum_{n=1}^\infty 1/p_n$
might converge.  Euler's result puts limitations on how quickly the
sequence of primes can grow.

Euler's proof uses the series~\eqref{eqn:zetaseries} for $s>1$ a real number.
He exploits the fact that the series grows without bound as $s\to 1^+$,
which follows from the fact that the harmonic series diverges.  But it is Riemann's name, not Euler's,
which is attached to the $\zeta$-function, and rightly so.
Riemann had the brilliant insight to consider $\zeta(s)$ for $s$ a complex number.

What does the series \eqref{eqn:zetaseries} mean when $s$ is not real?
Let $s=\sigma + i t$ be a complex variable, where $\sigma$ and $t$ are
real variables.  (The seemingly unusual choice of variable names 
has become standard when discussing the Riemann $\zeta$-function.)
The exponential of a complex number is defined by Euler's formula
$e^{i \theta} = \cos(\theta) + i \sin(\theta)$ for $\theta\in \mathbb R$.
Thus,
\begin{equation}
n^s = n^{\sigma + i t} = n^\sigma n^{it} = n^\sigma e^{i t \log(n)}
= n^\sigma (\cos(t \log(n)) + i \sin(t \log(n)) ).
\end{equation}

That formula explains what the terms in the series mean for complex~$s$, but
where does the series converge?
For any real number $\theta$ we have
\begin{equation}
| \cos(\theta) + i \sin(\theta) | = \sqrt{\cos^2(\theta) + \sin^2(\theta)} = \sqrt{1} = 1 .
\end{equation}
Thus,  $|n^s| = n^\sigma$.  Therefore, by the so-called\footnote{We write ``so-called'' because it was a terrible idea to name a result after the
variable in the expression, since the letter used for the variable is
an arbitrary choice.  The name ``zeta test'' is better.} ``$p$-series test'' from
calculus, which we will call the ``zeta test,''
\eqref{eqn:zetaseries} converges absolutely for $\sigma > 1$.

The Dirichlet series for $\zeta(s)$ only makes sense for $\sigma > 1$,
but there are other expressions for the $\zeta$-function, which  we will
not describe here, which allow us to define $\zeta(s)$ for all $s$ except
for $s=1$.  At $s=1$ the $\zeta$-function has a ``simple pole,'' which means
that $|\zeta(s)| \sim c/|s-1|$  as $s\to 1$.  In other words, the
$\zeta$-function ``blows up'' for $s$ near~$1$, but in the nicest possible way.

Riemann discovered that the zeros of the $\zeta$-function encode
secrets about the prime numbers.  The \emph{prime number theorem} is
the assertion that $p_n \sim n\log(n)$, meaning that 
$\lim_{n\to\infty} \frac{p_n}{n\log(n)} = 1$.  Here ``$\log$'' means
``natural logarithm.''  
Here is a silly proof of Euler's theorem that $\sum_{n=1}^\infty 1/p_n$ diverges:
combine the
prime number theorem, the integral test, and the limit
comparison test.
That proof is silly because the prime number theorem is a deep result,
and Euler's theorem can be proven with much more mundane input.

The prime number theorem was proven in 1896
by Hadamard and de la Vall\`ee Poussin, independently.
The key step in their proofs was
showing that $\zeta(s)$ has no zeros on the $\sigma=1$ line.
Riemann knew that the $\zeta$-function had no zeros for $\sigma > 1$.
The prime number theorem followed by turning that into a non-strict inequality.
That may seem like a minor improvement, but in fact it was one of the
major results of 19th century mathematics.

For the remainder of this paper we will set aside connections to the prime numbers and just
focus on analytic properties of the $\zeta$-function.

\subsection{Symmetries of the Riemann zeta-function}\label{sec:zeta}
The $\zeta$-function has a hidden symmetry which is revealed by the
Riemann $\xi$-function, defined by
\begin{equation}\label{eqn:xi}
\xi(s) = s (s-1) \pi^{-s/2} \Gamma(s/2) \zeta(s) .
\end{equation}
Here $\Gamma(s)$ is Euler's Gamma-function, defined by
\begin{equation}\label{eqn:Gamma}
\Gamma(s) = \int_0^\infty e^x x^s \frac{dx}{x} .
\end{equation}
The integral~\eqref{eqn:Gamma} converges for~$\sigma > 0$, and one can also
check that~$\Gamma(1)=1$.  Just as happened with
the Dirichlet series for $\zeta(s)$, we have an expression for $\Gamma(s)$
which is valid for some $s$, but we want to make use of the $\Gamma$-function for other
values of~$s$.  Euler figured out how to do that.
Integrating by parts (with $u=x^{s-1}$ and $dv = e^x dx$) we
find
\begin{equation}
\Gamma(s) = (s - 1) \Gamma(s-1) .
\end{equation}
which us usually written in the form (which is equivalent, because we are just
renaming variables) $\Gamma(s+1) = s \Gamma(s)$.  
If we rearrange that to get $\Gamma(s) = \Gamma(s+1)/s$, then we see that the
right side is known for $\sigma > -1$ except for $s=0$, and we can use that
as the definition for the left side.  Thus, we now know $\Gamma(s)$ for
$\sigma > -1$, except for $s=0$.  Repeating the same procedure we can define
$\Gamma(s)$ for $\sigma > -2$, except for $s=-1$ or $0$, then for
$\sigma > -3$, except for $s=-2$, $-1$, or $0$, and so on.

Euler's motivation was somewhat different:
he used the expression $\Gamma(s+1) = s \Gamma(s)$
to show that if $n$ is a positive integer, then $\Gamma(n) = (n-1)!$.
In other words, the $\Gamma$-function provides a way of extending the
factorial function to any positive real number.

Riemann showed that \eqref{eqn:xi} defines $\xi(s)$ for all $s\in\mathbb{C}$
and furthermore it satisfies the following surprising symmetry
which is called the \term{functional equation}:
\begin{equation}
\xi(s) = \xi(1-s) .
\end{equation}
In other words, the $\xi$-function is symmetric around the line $\sigma=\frac12$.
Combining the functional equation with the fact that $\xi(s)$ is real when
$s$ is real, we can conclude that
the zeros of the $\xi$-function either lie on the $\frac12$-line,
or occur in pairs symmetrically on either side of it.  The \term{Riemann hypothesis} is
the conjecture, made by Riemann in 1859, that all zeros of the $\xi$-function
actually lie on the $\sigma=\frac12$ line.  These days the line $\sigma=\frac12$ is known
as the \term{critical line}. 

Note that the zeros of $\xi(s)$ are also
zeros of $\zeta(s)$, but $\zeta(s)$ actually has some additional zeros which
are considered unimportant and go by the name~\term{trivial zeros}.  
The Riemann hypotheses is usually stated in the form ``the nontrivial zeros of
the Riemann zeta function lie on the critical line.''

We are one step away from connecting to our earlier discussion on real polynomials
with real zeros.  We can convert $\xi(s)$ to a function which is real on the
real axis by a simple change of variables, creating the Riemann $\Xi$-function:
\begin{equation}
\Xi(z) = \xi(\tfrac12 + i z) .
\end{equation}
Here $z=x + i y$ is a complex variable, and $x$ and $y$ are real.
We have: $\Xi(z)$ is real when $z$ is real, its zeros lie either on the real line
or in complex conjugate pairs, and the Riemann Hypotheses (RH) is the conjecture
that all zeros of~$\Xi(z)$ are real.  It is universally believed by all experts,
even those who are skeptical about RH, that all zeros of $\Xi(z)$ are simple.
Thus, we are in a situation parallel to our discussion about real polynomials
with distinct real zeros,
except that $\Xi(z)$ has infinitely many zeros.
Although we don't know that all the zeros of $\Xi$ are real,
we do know that all zeros either are real, or they lie close to the
real axis.  Combining the functional equation with the fact that $\zeta(s)$
is nonzero for $\sigma \ge 1$ shows that all zeros of $\Xi(z)$ lie
in the strip $-\frac12 < y < \frac12$.

One of the magical facts about polynomials is they can be written
as a product of linear factors,
each factor corresponds to a zero, and the zeros determine everything about
the polynomial except for a constant factor.  We will see that something
similar holds for certain functions that have infinitely many zeros.

\subsection{Another way to factor a polynomial}
Suppose $f(z) = a_0 + a_1 z + \cdots + a_n z^n$ is a polynomial
of degree~$n$, with roots $z_1$, \ldots, $z_n$.  If none of the
$z_j$ equal~$0$, then we can rearrange \eqref{eqn:polyfactored}  to get
\begin{align}\label{eqn:polyprod}
f(z) =\mathstrut & a_0 \left(1 - \frac{z}{z_1}\right)\cdots \left(1-\frac{z}{z_n}\right) \cr  
= &\mathstrut a_0 \prod_{j=1}^{n} \left(1 - \frac{z}{z_j}\right),
\end{align}
where in the second line we used product notation, which is analogous to
summation notation.
If $0$ is a root with multiplicity~$k$, and $z_{k+1}$, \ldots, $z_{n}$ are the
nonzero roots, then we have
\begin{equation}\label{eqn:polyprod0}
f(z) 
= \mathstrut a_k z^k \prod_{j=k+1}^{n} \left(1 - \frac{z}{z_j}\right).
\end{equation}

The product form~\eqref{eqn:polyprod} is useful because it can generalize
to the case of infinitely many zeros.  But before going into the details,
let's consider the concept of an infinite product.

\subsection{Infinite products}
Infinite sums
are familiar:  $\sum_{n=1}^\infty a_n$.  That sum might converge
or it might diverge.  It converges if the $a_n$ go to 0 fast enough.
Otherwise, it diverges.  An empty sum
is defined to equal~$0$.

Similarly, an infinite product, $\prod_{n=1}^\infty b_n$, might
converge, or it might diverge.  It converges if the $b_n$ go to 1
fast enough.  Otherwise, it diverges. 
The empty product is defined to equal~$1$.

The definition of convergent infinite product is analogous to
convergent infinite sum, with one subtle difference:
\begin{definition}[Convergence of an infinite product]
The infinite product
\begin{equation}
\prod_{n=1}^\infty b_n
\end{equation}
\term{converges} if there exists $n_0$ such that
$b_n$ is nonzero for $n \ge n_0$, 
and the limit
\begin{equation}\label{eqn:prodlim}
\lim_{N\to\infty} \prod_{n = n_0}^N b_n
\end{equation}
exists and is nonzero.
If $P$ is the limit in \eqref{eqn:prodlim}, then we say that the
infinite product converges to $b_1 \cdots b_{n_0 -1} P$.
\end{definition}

A possibly unexpected feature of the definition is that in a convergent
infinite product, only finitely many terms can be zero, and the product of
the nonzero terms also cannot equal zero.
That condition becomes natural
when relating infinite products to infinite sums.  Taking logarithms we have
\begin{equation}\label{eqn:logprod}
\log\left(\prod_{n = n_0}^N b_n \right) = \sum_{n = n_0}^N \log(b_n) .
\end{equation}
Since each $b_n$ in \eqref{eqn:logprod} is nonzero, and (for the infinite
product to have any hope of converging)  close to $1$ if
$n$ is sufficiently large, its logarithm is well-defined and close to~$0$. 
If the limit as $N\to\infty$ exists on the right side of \eqref{eqn:logprod},
then the limit inside the parentheses on the left side of \eqref{eqn:logprod}
exists \emph{and is nonzero}.  Thus, \emph{the infinite product
$\prod_{n=1}^\infty b_n$ converges if and only if the infinite sum
$\sum_{n=n_0}^\infty \log(b_n)$ converges}.

Before looking at the convergence of some infinite products,
let's recall the geometric series
\begin{equation}
\sum_{n=0}^\infty x^n = \frac{1}{1-x}
\end{equation}
which we can integrate to obtain
\begin{equation}
\sum_{n=0}^\infty\frac{x^{n+1}}{n+1} = -\log(1-x) .
\end{equation}
A useful way to rearrange that expression is
\begin{equation}
\log(1+x) = x - \tfrac12 x^2 + \text{smaller terms}.
\end{equation}
We will use this to determine whether certain infinite products converge.

Let $z$ be any number.  Does the product
\begin{equation}
\prod_{j=1}^\infty \left(1 - \frac{z}{j}\right)
\end{equation}
converge?  First: are the terms nonzero if $j$ is sufficiently large?
Yes, we just need $j > |z|$.  Second, taking logarithms we have
\begin{align}
\log\left( \prod_{j=j_0}^\infty \left(1 - \frac{z}{j}\right) \right)
=\mathstrut &
\sum_{j=j_0}^\infty \log \left(1 - \frac{z}{j}\right) \cr
=\mathstrut &
\sum_{j=j_0}^\infty \left( - \frac{z}{j} + \frac12 \frac{z^2}{j^2} + \text{smaller terms}\right) .
\end{align}
We recognize the harmonic series, which diverges, so the product diverges.

Does the product
\begin{equation}
\prod_{j=1}^\infty \left(1 - \frac{z}{j^2}\right)
\end{equation}
converge?  First: are the terms nonzero if $j$ is sufficiently large?
Yes, we just need $j > \sqrt{|z|}$.  Second, taking logarithms we have
\begin{align}
\log\left( \prod_{j=j_0}^\infty \left(1 - \frac{z}{j^2}\right) \right)
=\mathstrut &
\sum_{j=j_0}^\infty \log \left(1 - \frac{z}{j^2}\right) \cr
=\mathstrut &
\sum_{j=j_0}^\infty \left( - \frac{z}{j^2} + \tfrac12 \frac{z^2}{j^4} + \text{smaller terms}\right) .
\end{align}
By the zeta test, those series converge, therefore so does the infinite product.

Does the product
\begin{equation}\label{eqn:prodwithexp}
\prod_{j=1}^\infty \left(1 - \frac{z}{j}\right) e^{z/j}
\end{equation}
converge? As in the first example above, the terms are nonzero if $j > |z|$.
When we take logarithms, the summand is
\begin{align}\label{eqn:expfactor}
\log \left( \left(1 - \frac{z}{j}\right) e^{z/j} \right)
=\mathstrut &
\log  \left(1 - \frac{z}{j}\right) + \log(e^{z/j}) \cr
= \mathstrut & 
- \frac{z}{j} + \tfrac12 \frac{z^2}{j^2} + \text{smaller terms} + \frac{z}{j} \cr
= \mathstrut &  \tfrac12 \frac{z^2}{j^2} + \text{smaller terms} ,
\end{align}
so by the zeta test, the sum, and hence the product,
converges.

In the above discussion we have been somewhat imprecise in our use
of the approximation $\log(1 + x) = x - \tfrac12 x^2 + \cdots$.  The issue
which can arise is that just because a term is smaller, does not mean that an
infinite sum involving that term will be smaller.  For example,
consider
\begin{equation}
\sum_{n=1}^\infty \log\left(1 + \frac{(-1)^n}{\sqrt{\mathstrut n}}\right)
=\sum_{n=1}^\infty  \frac{(-1)^n}{\sqrt{\mathstrut n}} - \tfrac12 \sum_{n=1}^\infty \frac{1}{n} + \cdots .
\end{equation}
The first sum on the right converges (by the alternating series test).
The second sum has smaller terms (because $\frac{1}{\sqrt{n}} \ge \frac{1}{n}$),
but the second sum diverges because it is the harmonic series.
Thus, just because terms are smaller doesn't mean that we can ignore them.
But there is a case when we can:  when all the terms are positive.
That motivates this definition:
\begin{definition}\label{def:absconv}
The product
\begin{equation}
\prod_{n=1}^\infty (1 + c_n)
\end{equation}
\term{converges absolutely} if
\begin{equation}
\sum_{n=1}^\infty  c_n 
\end{equation}
converges absolutely.
\end{definition}
The ideas in the previous discussion prove that an absolutely convergent
product converges.

\section{Entire functions of finite order}
Polynomials are classified by their degree, or equivalently, the number of zeros
counting multiplicity.  There is another, also equivalent way, to classify polynomials:
how fast they grow.  If $f(z)=a_0 + \cdots + a_n z^n$ is a polynomial of degree~$n$,
and $|z|$ is very large, then $|f(z)|$ is approximately~$|a_n z^n|$.  
In other words, if $\mathcal N > \deg(f)$ and  $|z|$ is sufficiently large, then
\begin{equation}\label{eqn:degreebound}
\log(|f(z)|) < \mathcal N \log(|z|).
\end{equation}
We can use \eqref{eqn:degreebound} as an alternate definition of the degree of
a polynomial:  it is the infimum (greatest lower bound) of the numbers~$\mathcal N$
for which \eqref{eqn:degreebound} holds as $|z| \to \infty$.

The value of having multiple equivalent definitions is that some generalize
to other contexts, and some do not.  In this case, it is the rate of growth
which can be used to classify functions that have infinitely many zeros.

We will consider \term{entire} functions,
which means that the function has a Taylor series which converges for
all~$z\in \mathbb C$.  Examples are $\sin(z)$, $e^z$, polynomials, as well as
sums, products, and compositions of those functions.  Non-examples
are $\log(z)$, $1/z$, $\tan(z)$, and $\sqrt{z}$.
A good reference for the material in this section is~\cite{Lev}.

For polynomials there is a direct relationship between the growth of the
function and the number of zeros, but for entire functions the relationship
is not exact.  An entire function with many zeros must grow quickly, but
the function $e^z$ grows quickly even though it has no zeros.  It turns out
that is the complete story:  the growth of an entire function is
determined by its zeros, and by factors of the form~$e^{g(z)}$ where
$g$ is an entire function.  That motivates the definition of the
\term{order} of an entire function.

\begin{definition}  Suppose $f$ is an entire function, and suppose there
exists a real number $\mathcal N > 0$ such that if $z$ is sufficiently large,
\begin{equation}\label{eqn:orderN}
\log(|f(z)|) < |z|^{\mathcal N} .
\end{equation}
Then we say that $f$ \term{has finite order}.
If $f$ has finite order, then the \term{order} of~$f$ is the
infimum of the set of $\mathcal N$ such that~\eqref{eqn:orderN} holds
as $|z| \to \infty$.
\end{definition}

Combining~\eqref{eqn:degreebound} and \eqref{eqn:orderN}, and using the fact that
$\log(|z|)$ grows slower than any power of~$z$, we see that
polynomials have order~0.

The exponential function~$e^z$ has order~$1$,
and more generally if $g$ is a polynomial of degree~$N$ then
$e^{g(z)}$ has order~$N$.  Finite sums (or products) of functions of finite order
have finite order, and the order of the sum (or product) is at most
the order of the largest term (or factor).  The Riemann $\Xi$-function has
order~1, although that is not easy to deduce only from the information
we have provided in this paper.  The Euler $\Gamma$-function is not entire
because it is not defined at the negative integers, but $1/\Gamma(z)$ is an
entire function of order~1.

The beauty of entire functions of finite order is that they
can be written as an infinite product with a nice form.  This was
discovered by Hadamard, as part of his proof of the Prime Number Theorem.
We will state a special case which is sufficient for our needs.

\begin{theorem}[Hadamard]  Suppose $f$ is an entire function,
suppose~$0$ is a zero of $f$ of order~$m$, and
$z_1, z_2, \ldots$ are the non-zero zeros of~$f$.  If the order of $f$
is less than~1
then there exists $A\in \mathbb C$ such that
\begin{equation}\label{eqn:hadamard1}
f(z) = A z^m \prod_{n \ge 1} \left(1 - \frac{z}{z_n}\right) .
\end{equation}
If the order of $f$ is less than 2, then
 there exist numbers
$A, b \in \mathbb C$ such that
\begin{equation}\label{eqn:hadamard2}
f(z) = A e^{b z} z^m \prod_{n \ge 1} \left(1 - \frac{z}{z_n}\right) e^{z/z_n} .
\end{equation}
The products in \eqref{eqn:hadamard1} and \eqref{eqn:hadamard2}
converge absolutely for all~$z\in \mathbb C$.
\end{theorem}

By Definition~\ref{def:absconv},
the sum $\sum 1/|z_n|$ in \eqref{eqn:hadamard1} converges,
and the same ideas we used when discussing \eqref{eqn:prodwithexp} show that
in \eqref{eqn:hadamard2} the sum 
$\sum 1/|z_n|^2$ converges.
In general, for a function of order~$\rho$
with non-zero zeros $z_1, z_2, \ldots$, the series $\sum 1/|z_n|^p$ converges
if~$p > \rho$.  Thus there is a close relationship between the order of a function
and the distribution of its zeros, with low order functions having zeros
distributed more sparsely.

Now we have the ingredients to extend our results on polynomials to
functions with infinitely many zeros.

\subsection{Zeros of derivatives}
Entire functions of order less than 2 are nice because the analogue of 
Theorem~\ref{thm:easydirection} holds for them:
\begin{theorem}\label{thm:easyorder1} Suppose $f(z)$ is an entire function of 
order less than~2, which is real
on the real axis and has only distinct real zeros.  Then all zeros of $f'$ are real,
and furthermore all local maxima of $f'$ are positive and all local
minima of $f'$ are negative.
\end{theorem}

Since the derivative of an entire function of order $\rho$ is entire
of order~$\rho$, the conclusion of Theorem~\ref{thm:easyorder1} 
applies to every derivative~$f^{(n)}$ and not only~$f'$.


The proof of Theorem~\ref{thm:easydirection} used Rolle's theorem and the fact
that the derivative of a polynomial has one fewer zero than the original
polynomial.  
However, it makes
no sense to talk about ``one fewer'' zero of the derivative when there
are infinitely many zeros.  Instead,
we will show that $f'$ has exactly one real zero between each pair of
consecutive real zeros of $f$, and no other zeros. 

\begin{proof}
Suppose $f$ satisfies the conditions in the Theorem~\ref{thm:easyorder1}.
Since $f$ is an entire function of order less than~$2$, it has the 
form~\eqref{eqn:hadamard2}.  Furthermore, $A$, $b$, and all $z_j$ are real.
Since $f'(z)/f(z) = \frac{d}{dz}\log(f(z))$, 
we have
\begin{equation}\label{eqn:logderiv}
\frac{f'}{f}(z) = B + \sum_{j=1}^\infty \frac{1}{z - z_j} + \frac{1}{z_j} .
\end{equation}
Note that the series in~\eqref{eqn:logderiv} is defined whenever $z$ is not
a zero of $f$, and the sum converges absolutely wherever it is defined.

Suppose $z_j$ and $z_{j*}$ are consecutive zeros of $f$, and
$z_j < z < z_{j*}$.  We wish to show that $f'(z)$ has only
one zero in that interval.  Since $f(z)$ has constant sign on that interval,
this is equivalent to $f'(z)/f(z)$ having exactly one zero on that
interval.  That function has at least one zero on that interval because by
\eqref{eqn:logderiv}, that function is large and positive if $z$ is slightly
larger than $z_j$, and it is large and negative if $z$ is slightly
smaller than $z_{j*}$.  

To show that $f'(z)$, equivalently $f'(z)/f(z)$, has only one zero in
that interval, it is sufficient to show that $f'(z)/f(z)$ is decreasing
on that interval.  But that follows from differentiating~\eqref{eqn:logderiv}:
\begin{equation}
\left(\frac{f'}{f}(z)\right)' = \sum_{j=1}^\infty -\frac{1}{(z - z_j)^2},
\end{equation}
which is negative for real $z$ whenever it is defined.
In particular it is negative on the interval $(z_j,z_{j*})$.

We have shown that $f'(z)$ has only one zero in each interval $(z_j, z_{j*})$,
which proves that all maxima of $f$ are positive  and all minima of
$f$ are negative.  Now we must prove that $f'$ has no non-real zeros.
In \eqref{eqn:logderiv} let $z=x + i y$, and then take the imaginary
part of that expression, making use of the fact that $B$ and $z_j$ are real:
\begin{equation}
\Im\left(\frac{f'}{f}(x + i y)\right) =
\sum_{j=1}^\infty - \frac{y}{(x-z_j)^2 + y^2} .
\end{equation}
That is nonzero if $y\not=0$, because every term is nonzero and has
the same sign, therefore $f'(z)/f(z)$ cannot be $0$ if
$z$ is not real, because its imaginary part is not zero.
\end{proof}

The above proof is not quite right.  The first two paragraphs implicitly
assumed that the zeros of $f$ range from $-\infty$ to $\infty$, but that
might not be the case:  there might be a largest, or a smallest, zero,
or there might only be finitely many zeros.  We leave it as an exercise to
work out the details for those cases.  Illustrative examples to consider
are $e^{z} p(z)$ where $p$ is a real polynomial.

\subsection{The converse}
Does the analogue of Theorem~\ref{thm:harddirection} hold for entire
functions of order less than~2?
Definitely not.  The function $e^z - 1$ real on the real line and is entire
of order~$1$.  It has zeros at $2\pi i n$ for all $n\in \mathbb Z$,
so it has infinitely many complex zeros.
But it has no real maxima or minima.  Its first derivative, and
every higher derivative, is $e^z$, which also has no real maxima or minima.

All is not lost:  our motivation was the Riemann $\Xi$-function,
which has order~1 and is real on the real axis.  It has one more
key property:  its complex zeros, if they exist, must have imaginary
part between $-\frac12$ and~$\frac12$, which follows from the fact that
$\zeta(s)$ has no zeros in~$\sigma \ge 1$.  That is, all zeros lie in a
strip around the real axis.   The counterexample $e^z - 1$ has zeros
arbitrarily far from the real axis.  Maybe the analogue of
Theorem~\ref{thm:harddirection} holds for functions with all zeros
near the real axis?

Y-.O.~Kim~\cite{Kim} proved a version with that extra assumption,
plus an additional extra assumption.

\begin{theorem}\label{thm:kimhalf}  Suppose $f(z)$ is an entire function of order~1
which is real on the real axis and has all zeros in the strip
$-A < y < A$ for some $A > 0$.  Let $w_1, w_2, \ldots$ be the
non-real zeros of $f$ and further suppose
\begin{equation}\label{eqn:sumthird}
\sum_j \frac{1}{|w_j|^p} \text{ converges for some $p < \frac12$. }
\end{equation}
If $f$ and all of its derivatives have every maximum positive and
every minimum negative, then $f$ has only real zeros, and all zeros
are distinct.
\end{theorem}

In other words, if the zeros are in a strip and there are not too many
complex zeros, then the analogue of Theorem~\ref{thm:harddirection}
is true.  



Might Theorem~\ref{thm:kimhalf} be true without the extra condition on
the non-real zeros?  Or maybe if there are many non-real zeros, they can
somehow conspire to stay away from the real axis as the function is differentiated?
Those are research questions.

One could apply Theorem~\ref{thm:kimhalf} in its current form, of one could show
that the Riemann $\Xi$-function had sufficiently few non-real zeros.  This is
far beyond the realm of possibility in the foreseeable future.  At present
it has not even been shown that at least half of the zeros of the $\Xi$-function
are real.

At the start of this section we asked whether there might be an equivalence
to the Riemann Hypothesis that only involved the maxima and minima of
the $\Xi$-function and its derivatives.  Such an equivalence cannot follow
from a general statement about real entire functions of order~1,
as illustrated by the example
and $\cos(a z)(z^2 + \frac14)$.
That function has zeros at $\pm i/2$ and all other zeros are real.
As $a$ increases, the real zeros move closer together but the pair of
complex zeros do not move.
If $a$ is small,
the function has a positive local minimum, but if $a$ is large
enough, it does not.  As $a$ increases, it takes more and more derivatives
before a positive local minimum or negative local maximum appears.
That is a another version of a lesson from Figure~\ref{fig:motivatingplots}:
funny business comes from complex zeros, complex zeros close to the
real axis cause more funny business, and closeness to the real axis
should be measured relative to the gaps between nearby zeros.



\end{document}